\documentclass{amsart}

\usepackage{amsthm}
\usepackage[leqno]{amsmath}
\usepackage[all]{xy} \SelectTips{eu}{}
\usepackage{latexsym,amsfonts,amssymb}
\usepackage{appendix}
\usepackage{hyperref}

\newtheorem{intthm}{Theorem}[]

\newcommand{\numberseries}{\bfseries}   
\newlength{\thmtopspace}                
\newlength{\thmbotspace}                
\newlength{\thmheadspace}               
\newlength{\thmindent}                  

\newtheoremstyle{bfupright head,slanted body}
                {\thmtopspace}{\thmbotspace}
                {\slshape}{\thmindent}{\bfseries}{.}{\thmheadspace}
                {{\numberseries \thmnumber{#2\;}}\thmnote{#3}}

\newtheoremstyle{bfupright head,upright body}
                {\thmtopspace}{\thmbotspace}
                {\upshape}{\thmindent}{\bfseries}{.}{\thmheadspace}
                {{\numberseries \thmnumber{#2\;}}\thmnote{#3}}

\newtheoremstyle{fixed bf head,slanted body}
                {\thmtopspace}{\thmbotspace}{\slshape}
                {\thmindent}{\bfseries}{.}{\thmheadspace}
                {{\numberseries \thmnumber{#2\;}}\thmname{#1}\thmnote{ (#3)}}

\newtheoremstyle{fixed bf head,upright body}
                {\thmtopspace}{\thmbotspace}{\upshape}
                {\thmindent}{\bfseries}{.}{\thmheadspace}
                {{\numberseries \thmnumber{#2\;}}\thmname{#1}\thmnote{ (#3)}}

\newtheoremstyle{numbered paragraph}
                {\thmtopspace}{\thmbotspace}{\upshape}
                {\thmindent}{\upshape}{}{\thmheadspace}
                {{\numberseries \thmnumber{#2.}}}

\theoremstyle{bfupright head,slanted body}
\newtheorem{res}{}[section]             \newtheorem*{res*}{}

\theoremstyle{bfupright head,upright body}
               \newtheorem*{bfhpg*}{}

\theoremstyle{fixed bf head,slanted body}
\newtheorem{thm}[res]{Theorem}          \newtheorem*{thm*}{Theorem}
\newtheorem{prp}[res]{Proposition}      \newtheorem*{prp*}{Proposition}
\newtheorem{cor}[res]{Corollary}        \newtheorem*{cor*}{Corollary}
\newtheorem{lem}[res]{Lemma}            \newtheorem*{lem*}{Lemma}

\theoremstyle{fixed bf head,upright body}
\newtheorem{dfn}[res]{Definition}       \newtheorem*{dfn*}{Definition}
           \newtheorem*{rmk*}{Remark}
          \newtheorem*{exa*}{Example}

\theoremstyle{numbered paragraph}
\newtheorem{ipg}[res]{}

\newlength{\thmlistleft}        
\newlength{\thmlistright}       
\newlength{\thmlistpartopsep}   
\newlength{\thmlisttopsep}      
\newlength{\thmlistparsep}      
\newlength{\thmlistitemsep}     

\newcounter{eqc}
\newenvironment{eqc}{\begin{list}{\upshape (\textit{\roman{eqc}})}%
    {\usecounter{eqc}%
      \setlength{\leftmargin}{\thmlistleft}%
      \setlength{\labelwidth}{\thmlistleft}%
      \setlength{\rightmargin}{\thmlistright}%
      \setlength{\partopsep}{\thmlistpartopsep}%
      \setlength{\topsep}{\thmlisttopsep}%
      \setlength{\parsep}{\thmlistparsep}%
      \setlength{\itemsep}{\thmlistitemsep}}}%
  {\end{list}}%

\newcommand{\eqclbl}[1]{{\upshape(\textit{#1})}}

\newcounter{rqm}
  {\end{list}}%

  \newcommand{\proofofimp}[3][:]{\mbox{\eqclbl{#2}$\!\implies\!$\eqclbl{#3}#1}}
\newenvironment{prf*}[1][Proof]{%
  \begin{proof}[\bf #1]
    \setcounter{equation}{0}
    }
  {\end{proof}
}

\newcommand{\pgref}[1]{\ref{#1}}

\renewcommand{\eqref}[1]{(\pgref{eq:#1})}

\def\urltilda{\kern -.15em\lower .7ex\hbox{\~{}}\kern .04em}

\newcommand{\ZZ}{\mathbb{Z}}
\newcommand{\deq}{\:=\:}
\newcommand{\dis}{\:\is\:}
\DeclareMathOperator*{\colim}{colim}

\newcommand{\is}{\cong}
\newcommand{\qis}{\simeq}

\newcommand{\xra}[2][]{\xrightarrow[#1]{\;#2\;}}
\newcommand{\qra}{\xra{\smash{\qis}}}
\newcommand{\Rop}{R^\circ}
\newcommand{\Sop}{S^\circ}
\newcommand{\del}{\partial}

\newcommand{\dmapdef}[4][\lra]{\nobreak{#2\colon #3\:#1\:#4}}
\newcommand{\Ker}[1]{\nobreak{\operatorname{Ker}#1}}
\renewcommand{\Im}[1]{\nobreak{\operatorname{Im}#1}}
\newcommand{\Coker}[1]{\nobreak{\operatorname{Coker}#1}}
\renewcommand{\H}[2][]{\operatorname{H}_{#1}(#2)}

\newcommand{\Hom}[3][R]{\operatorname{Hom}_{#1}(#2,#3)}
\newcommand{\FHom}[3][R]{\operatorname{FHom}_{#1}(#2,#3)}
\newcommand{\Ext}[4][R]{\operatorname{Ext}_{#1}^{#2}(#3,#4)}
\newcommand{\tp}[3][R]{\nobreak{#2\otimes_{#1}#3}}

\makeatletter
\def\@nobreak@#1{\mathchoice%
  {\nobreakdef@\displaystyle\f@size{#1}}%
  {\nobreakdef@\nobreakstyle\tf@size{\firstchoice@false #1}}%
  {\nobreakdef@\nobreakstyle\sf@size{\firstchoice@false #1}}%
  {\nobreakdef@\nobreakstyle\ssf@size{\firstchoice@false #1}}%
  \check@mathfonts}%
\def\nobreakdef@#1#2#3{\hbox{{%
                    \everymath{#1}%
                    \let\f@size#2\selectfont%
                    #3}}}%
\makeatother

\def\widebardisplay#1{%
  \setbox0=\hbox{$\displaystyle #1$}
  \dimen0=\wd0%
  \advance\dimen0 by -2pt
  \vbox{%
    \nointerlineskip%
    \moveright 1pt 
    \vbox{\hrule width \dimen0}%
    \nointerlineskip%
    \kern 2pt
    \box0%
    }%
  }

\def\widebartext#1{%
  \setbox0=\hbox{$#1$}
  \dimen0=\wd0%
  \advance\dimen0 by -2pt
  \vbox{%
    \nointerlineskip%
    \moveright 1pt 
    \vbox{\hrule width \dimen0}%
    \nointerlineskip%
    \kern 1.6pt
    \box0%
    }%
  }

\def\widebarscript#1{%
  \setbox0=\hbox{$\scriptstyle #1$}
  \dimen0=\wd0%
  \advance\dimen0 by -3pt
  \vbox{%
    \nointerlineskip%
    \moveright 1.5pt 
    \vbox{\hrule width \dimen0}%
    \nointerlineskip%
    \kern .8pt
    \box0%
    }%
  }

\def\widebarscriptscript#1{%
  \setbox0=\hbox{$\scriptscriptstyle #1$}
  \dimen0=\wd0%
  \advance\dimen0 by -2pt
  \vbox{%
    \nointerlineskip%
    \moveright 1pt 
    \vbox{\hrule width \dimen0}%
    \nointerlineskip%
    \kern .6pt
    \box0%
    }%
  }

\def\widebar#1{\mathchoice%
  {\widebardisplay{#1}}%
  {\widebartext{#1}}%
  {\widebarscript{#1}}%
  {\widebarscriptscript{#1}}%
  }

\newcommand{\btp}[3][R]{\nobreak{#2\mathbin{\widebar{\otimes}}_{#1}#3}}

\newcommand{\ttp}[3][R]{\nobreak{#2\mathbin{\widetilde{\otimes}}_{#1}#3}}

\newcommand{\Stor}[4][R]{\smash{\operatorname{\widetilde{Tor}}}_%
  {#2}^{{#1}^{\phantom{|\mspace{-6mu}}}}(#3,#4)}

\newcommand{\sExt}[4][R]{\widetilde{\operatorname{Ext}}_{#1}^{#2}(#3,#4)}
\newcommand{\sHom}[4][R]{\widetilde{\operatorname{Hom}}_{#1}(#2,#3)}

\newcommand{\susp}{\mathsf{\Sigma}}
\newcommand{\SHom}[3][R]{\widetilde{\operatorname{Hom}}_{#1}(#2,#3)}
\newcommand{\bHom}[3][R]{\widebar{\operatorname{Hom}}_{#1}(#2,#3)}

\DeclareFontFamily{U}{mathx}{\hyphenchar\font45}
\DeclareFontShape{U}{mathx}{m}{n}{ <5> <6> <7> <8> <9> <10> <10.95>
  <12> <14.4> <17.28> <20.74> <24.88> mathx10 }{}
\DeclareSymbolFont{mathx}{U}{mathx}{m}{n}
\DeclareFontSubstitution{U}{mathx}{m}{n}
\DeclareMathAccent{\widecheck}{0}{mathx}{"71}
\DeclareMathAccent{\wideparen}{0}{mathx}{"75}

\newcommand{\pd}{\operatorname{pd}}
\newcommand{\id}{\operatorname{id}}
\newcommand{\fd}{\operatorname{fd}}

\begin{document}
\title[Stable homology with respect to a semidualizing module]{Vanishing of stable homology with respect to a semidualizing module}

\author[L. Liang]{Li Liang}

\address{Li Liang, School of Mathematics and Physics, Lanzhou Jiaotong University, Lanzhou 730070, China.}

\email{lliangnju@gmail.com}

\thanks{This research was partly supported by NSF of China (Grant Nos. 11761045, 11301240 and 11561039), NSF of Gansu Province (Grant No. 1506RJZA075) and SRF for ROCS, SEM}

\date{\today}



\begin{abstract}
We investigate stable homology of modules over a commutative noetherian ring $R$ with respect to a semidualzing module $C$, and give some vanishing results that improve/extend the known results. As a consequence, we show that the balance of the theory forces $C$ to be trivial and $R$ to be Gorenstein.
\end{abstract}

\maketitle

\thispagestyle{empty} \allowdisplaybreaks

\section{Introduction}

\noindent
Stable homology, as a broad generalization of Tate homology to the realm of associative rings, was introduced by Vogel and Goichot \cite{Go}, and further studied by Celikbas, Christensen, Liang and Piepmeyer \cite{CCLP01,CCLP02}, and Emmanouil and Manousaki \cite{EM}. In their paper \cite{CCLP01}, it is shown that the vanishing of stable homology over commutative noetherian local rings can detect modules of finite projective (injective) dimension, even of finite Gorenstein dimension, which lead to some characterizations of classical rings such as Gorenstein rings, the original domain of Tate homology. In \cite{EM}, Emmanouil and Manousaki further investigate stable homology of modules, and give some vanishing results that improve results in \cite{CCLP01} by relaxing the conditions on rings and modules.

The study of semidualizing modules was initiated independently by Foxby \cite{F}, Golod \cite{Gol} and Vasconcelos \cite{V}. Over a commutative noetherian ring $R$, a finitely generated $R$-module $C$ is semidualizing if $\Hom{C}{C} \is R$ and $\Ext{i}{C}{C}=0$ for all $i\geq 1$. Examples include finitely generated projective $R$-modules of rank $1$. Modules of finite homological dimension with respect to a semidualizing module have been studied in numerous papers. For example, Takahashi and White \cite{TW}, and Salimi, Sather-Wagstaff, Tavasoli and Yassemi \cite{SSWTY} give some characterizations for such modules in terms of the vanishing of relative (co)homology. In this paper, we show that the vanishing of stable homology can also detect modules of finite homological dimension with respect to a semidualizing module. Our main results are two theorems as shown below; see Theorems \ref{cinj dim-Tor} and \ref{cflat dim-Tor}.

\begin{intthm}\label{cflat}
Let $R$ be a commutative noetherian ring and let $C$ be a semidualizing $R$-module. For an $R$-module $M$, the following conditions are equivalent.
  \begin{eqc}
  \item $\mathcal{F}_C{\text-}\pd_{R}M<\infty$.
  \item $\Stor[\mathcal{P}_C\mathcal{I}_C]{n}{M}{-}=0$ for each $n\in\ZZ$.
  \item $\Stor[\mathcal{P}_C\mathcal{I}_C]{n}{M}{-}=0$ for some $n\geq0$.
  \end{eqc}
\end{intthm}

\begin{intthm}\label{cinj}
Let $R$ be a commutative noetherian ring and let $C$ be a semidualizing $R$-module. For an $R$-module $N$, the following conditions are equivalent.
  \begin{eqc}
  \item $\mathcal{I}_C{\text-}\id_{R}N<\infty$.
  \item $\Stor[\mathcal{P}_C\mathcal{I}_C]{n}{-}{N}=0$ for each $n\in\ZZ$.
  \item $\Stor[\mathcal{P}_C\mathcal{I}_C]{n}{-}{N}=0$ for some $n<0$.
  \end{eqc}
\end{intthm}

The above two results improve the right and left vanishing results in the introduction of \cite{CCLP01}. Here the notation $\mathcal{F}_C{\text-}\pd_{R}M$, $\mathcal{I}_C{\text-}\id_{R}N$ and $\Stor[\mathcal{P}_C\mathcal{I}_C]{n}{-}{-}$ can be found in \ref{resolution} and \ref{stable homology}. As a consequence, we show that the isomorphisms $\Stor[\mathcal{P}_C\mathcal{I}_C]{\ast}{M}{N} \is \Stor[\mathcal{P}_C\mathcal{I}_C]{\ast}{N}{M}$ for all $R$-modules $M$ and $N$ force $C$ to be trivial and $R$ to be a Gorenstein ring; see Corollary \ref{c trival}.

We prove these results using the next characterization of stable (unbounded) tensor product inspired by the work of Emmanouil and Manousaki \cite{EM}; see Theorem \ref{stensor}.

\begin{intthm}
Let $X$ be a complex of $\Rop$-modules and $Y$ a bounded above complex of $R$-modules with $\sup\{i\in\ZZ|Y_{i}\neq 0\}=k$. Then there are isomorphisms of complexes of $\mathbb{Z}$-modules
$$\btp{X}{Y} \is {\lim}_{i\in\mathbb{N}}((\tp{X}{Y})/(\tp{X}{Y_{\leq k-i}})),$$
and
$$\ttp{X}{Y} \is {\lim}^{1}_{i\in\mathbb{N}}(\tp{X}{Y_{\leq k-i}}).$$
\end{intthm}

One refers to \ref{sten def} for the definitions of $\btp{X}{Y}$ and $\ttp{X}{Y}$, and ${\lim}^{1}$ is the right derived functor of the limit lim; see \ref{limits}.

\section{Preliminaries}

\noindent
We begin with some notation and terminology for use throughout this paper.

\begin{ipg}
Throughout this work, all rings are assumed to be associative rings. Let $R$ be a ring; by an $R$-module we mean a left $R$-module, and we refer to right $R$-modules as modules over the opposite ring $\Rop$. We denote by $\mathcal{P}$ (resp., $\mathcal{F}$, $\mathcal{I}$) the class of projective $R$-modules (resp., flat $R$-modules, injective $R$-modules).

By an $R$-complex we mean a complex of $R$-modules. We frequently (and without warning) identify $R$-modules with $R$-complexes concentrated in degree $0$. For an $R$-complex $X$, we set $\sup X=\sup\{i\in\ZZ|X_{i}\neq 0\}$ and $\inf X=\inf\{i\in\ZZ|X_{i}\neq 0\}$. An $R$-complex $X$ is \emph{bounded above} if $\sup X<\infty$, and it is \emph{bounded below} if $\inf X>-\infty$. An $R$-complex $X$ is \emph{bounded} if it is both bounded above and bounded below. The $n$th \emph{homology} of $X$ is denoted by $\H[n]{X}$.
For each $k\in\mathbb{Z}$, $\susp^{k}X$ denotes the complex with the degree-$n$ term $(\susp^{k}X)_{n}=X_{n-k}$ and whose boundary operators are $(-1)^{k}\del_{n-k}^{X}$. We set $\susp M=\susp^{1}M$.

If $X$ and $Y$ are both $R$-complexes, then by a \emph{morphism} $\alpha:\xymatrix@C=0.5cm{X \ar[r] & Y}$ we mean a sequence $\alpha_{n}:\xymatrix@C=0.5cm{X_{n} \ar[r] & Y_{n}}$ such that $\alpha_{n-1}\del^{X}_{n}=\del^{Y}_{n}\alpha_{n}$ for each $n\in\mathbb{Z}$. A \emph{quasi-isomorphism}, indicated by the symbol ``$\simeq$", is a morphism of complexes that induces an isomorphism in homology.
\end{ipg}

\begin{ipg} \label{resolution}
Let $\mathcal{X}$ be a class of $R$-modules. Following Enochs and Jenda \cite{EJ}, an $\mathcal{X}$-precover of an $R$-module $M$ is a homomorphism $X\to M$ with $X\in\mathcal{X}$ such that the homomorphism $\Hom{X'}{X} \to \Hom{X'}{M}$ is surjective for each $X'\in\mathcal{X}$. $\mathcal{X}$ is called a precovering class if each $R$-module has a $\mathcal{X}$-precover.

For a precovering class $\mathcal{X}$, there is a complex
$$X^{+}\equiv\xymatrix@C=0.5cm{\cdots \ar[r] & X_{1} \ar[r] & X_{0} \ar[r] & M \ar[r] & 0 }$$
with each $X_{i}$ in $\mathcal{X}$, such that $\Hom{X'}{X^{+}}$ is exact for each $X'\in\mathcal{X}$. The truncated complex $X\equiv\xymatrix@C=0.5cm{\cdots \ar[r] & X_{1} \ar[r] & X_{0} \ar[r] & 0 }$ is called a \emph{proper $\mathcal{X}$-resolution} of $M$, which is always denoted by $X \to M$. If $\mathcal{X}$ contains all projective $R$-modules, then the complex $X^{+}$ is exact. In this case, we always denote by $X \qra M$ the proper $\mathcal{X}$-resolution of $M$.

The \emph{$\mathcal{X}$-projective dimension} of $M$ is the quantity:
$$\mathcal{X}{\text-}\pd_{R}{M}=\inf\{\sup X\ |\ X\to M\ {\rm is\ a\ proper}\ \mathcal{X}{\text-}{\rm resolution\ of}\ M\}.$$

We define \emph{preenveloping classes $\mathcal{Y}$}, \emph{proper $\mathcal{Y}$-coresolutions} and \emph{$\mathcal{Y}$-injective dimension}, $\mathcal{Y}{\text-}\id_{R}M$, of $M$ dually.

When $\mathcal{X}$ is the class of projective (resp, flat) $R$-modules, $\mathcal{X}{\text-}\pd_{R}{M}$ is the classical projective (resp. flat) dimension; we refer the reader to \cite[Remark 2.6]{SSWTY} for the flat case. Also when $\mathcal{Y}$ is the class of injective $R$-modules, $\mathcal{Y}{\text-}\id_{R}M$ is the classical injective dimension.
\end{ipg}

\section{A characterization of stable (unbounded) tensor product} 

\noindent
We start by recalling the definition of stable (unbounded) tensor product.

\begin{ipg}\label{sten def}
Let $X$ be an $\Rop$-complex and $Y$ an $R$-complex. The \emph{tensor product} $\tp{X}{Y}$ is the $\mathbb{Z}$-complex with degree-$n$ term
$$(X\otimes_{R}Y)_{n}=\coprod_{i\in\ZZ}(X_{i}\otimes_{R}Y_{n-i})$$
and differential given by $\del^{\tp{X}{Y}}(x\otimes y)=\del^{X}(x)\otimes
y+(-1)^{|x|}x\otimes\del^{Y}(y)$. Following \cite{CCLP01,Go}, the \emph{unbounded tensor product} $\btp{X}{Y}$ is the $\mathbb{Z}$-complex with degree-$n$ term
$$(X\otimes_{R}Y)_{n}=\prod_{i\in\ZZ}(X_{i}\otimes_{R}Y_{n-i})$$
and differential defined as above. $\tp{X}{Y}$ is a subcomplex of $\btp{X}{Y}$, so we let $\ttp{X}{Y}$ denote the quotient complex $(\btp{X}{Y})/(\tp{X}{Y})$, which is called the \emph{stable tensor product}.
\end{ipg}

We notice that if $X$ or $Y$ is bounded, or if both of them are bounded on the same
side (above or below), then the unbounded tensor product coincides with the
tensor product, and so the stable tensor product $\ttp{X}{Y}$ is zero.

\begin{ipg}\label{limits}
Let ${\{\nu^{uv}: X^{v} \to X^{u}}\}_{u\leq v}$ be an $\mathbb{N}$-inverse system of $R$-complexes. For the morphism $1-\nu: \prod_{i\in\mathbb{N}}X^{i} \to \prod_{i\in\mathbb{N}}X^{i}$ given by $(1-\nu)_{k}(x_{i})_{i\in\mathbb{N}}=(x_{i}-\nu_{k}^{i,i+1}(x_{i+1}))_{i\in\mathbb{N}}$ for each $k\in\ZZ$, where $(x_{i})_{i\in\mathbb{N}}\in\prod_{i\in\mathbb{N}}X_{k}^{i}$, it is well known that $\Ker(1-\nu)=\lim_{i\in\mathbb{N}}X^{i}$ and $\Coker(1-\nu)=\lim^{1}_{i\in\mathbb{N}}X^{i}$. Here $\lim^{1}$ is the right derived functor of the limit lim; see e.g. Emmanouil \cite{E}, Roos \cite{R} and Yeh \cite{Y} for more details. That is, there is an exact sequence of $R$-complexes
$$0 \to {\lim}_{i\in\mathbb{N}}X^{i} \to \prod_{i\in\mathbb{N}}X^{i} \to \prod_{i\in\mathbb{N}}X^{i} \to {\lim}^1_{i\in\mathbb{N}}X^{i} \to 0.$$
\end{ipg}

Let $X$ be an $R$-complex, and $X=X^{0}\supseteq X^1\supseteq \cdots$ a filtration. Then the embeddings $\varepsilon^{i}: X^{i} \to X^{i-1}$ and the morphisms $\pi^{i}: X/X^i \to X/{X^{i-1}}$ determine the $\mathbb{N}$-inverse systems
$${\{\varepsilon^{uv}: X^{v} \to X^{u}}\}_{u\leq v}\ {\rm and} \ {\{\pi^{uv}: X/X^{v} \to X/X^{u}}\}_{u\leq v}\,$$
respectively. For these systems, we have the following result.

\begin{lem}\label{exact sequence}
Let $X$ be an $R$-complex, and $X=X^{0}\supseteq X^1\supseteq \cdots$ a filtration. Then $\lim^{1}_{i\in\mathbb{N}}X/X^{i}=0$, and there exists an exact sequence
$$0 \to {\lim}_{i\in\mathbb{N}}X^{i} \to X \to {\lim}_{i\in\mathbb{N}}X/X^{i} \to {\lim}^{1}_{i\in\mathbb{N}}X^{i} \to 0.$$
\end{lem}
\begin{prf*}
Consider the following commutative diagram with exact rows
$$\xymatrix{
  0 \ar[r] & \prod_{i\in\mathbb{N}}X^{i} \ar[d]^{1-\varepsilon} \ar[r] & \prod_{i\in\mathbb{N}}X \ar[d]^{1-\id} \ar[r] & \prod_{i\in\mathbb{N}}X/X^{i} \ar[d]^{1-\pi} \ar[r] & 0 \\
  0 \ar[r] & \prod_{i\in\mathbb{N}}X^{i} \ar[r] & \prod_{i\in\mathbb{N}}X \ar[r] & \prod_{i\in\mathbb{N}}X/X^{i} \ar[r] & 0.}$$
We notice that the constant $\mathbb{N}$-inverse system $\{X\}$ has $\lim_{i\in\mathbb{N}}X=X$ and $\lim^1_{i\in\mathbb{N}}X=0$ since $1-\id$ is surjective. Then by \ref{limits} and the snake lemma, one gets the desired results.
\end{prf*}

\begin{ipg}\label{filtration}
Let $X$ be an $\Rop$-complex and $Y$ an $R$-complex. Fix $k\in\mathbb{Z}$, the
filtration
$$Y_{\leq k}\supseteq Y_{\leq k-1}\supseteq Y_{\leq k-2}\supseteq \cdots$$
induces a filtration
$$\tp{X}{Y_{\leq k}}\supseteq \tp{X}{Y_{\leq k-1}}\supseteq \tp{X}{Y_{\leq k-2}}\supseteq \cdots.$$
Thus we have two $\mathbb{N}$-inverse systems $\{\varepsilon^{uv}: \tp{X}{Y_{\leq k-v}} \to \tp{X}{Y_{\leq k-u}}\}_{u\leq v}$ and $\{\pi^{uv}: (\tp{X}{Y_{\leq k}})/(\tp{X}{Y_{\leq k-v}}) \to (\tp{X}{Y_{\leq k}})/(\tp{X}{Y_{\leq k-u}})\}_{u\leq v}$.
\end{ipg}

\begin{thm}\label{stensor}
Let $X$ be an $\Rop$-complex and $Y$ a bounded above $R$-complex with $\sup Y=k$. Then there are isomorphisms of $\mathbb{Z}$-complex
$$\btp{X}{Y} \is {\lim}_{i\in\mathbb{N}}((\tp{X}{Y})/(\tp{X}{Y_{\leq k-i}})),$$
and
$$\ttp{X}{Y} \is {\lim}^{1}_{i\in\mathbb{N}}(\tp{X}{Y_{\leq k-i}}).$$
\end{thm}
\begin{prf*}
We first prove the case where $k=0$. In this case, $Y=Y_{\leq 0}$. For each $n\in\mathbb{Z}$,
$$(\tp{X}{Y_{\leq 0}})_{n}=\coprod_{p\in\mathbb{Z}}(\tp{X_{n+p}}{(Y_{\leq 0})_{-p}})=\coprod_{p\geq 0}(\tp{X_{n+p}}{(Y_{\leq 0})_{-p}}),$$
and for each $i\geq 1$,
$$(\tp{X}{Y_{\leq -i}})_{n}=\coprod_{p\in\mathbb{Z}}(\tp{X_{n+p}}{(Y_{\leq -i})_{-p}})=\coprod_{p\geq i}(\tp{X_{n+p}}{(Y_{\leq 0})_{-p}}).$$
Thus one gets
$$((\tp{X}{Y_{\leq 0}})/(\tp{X}{Y_{\leq -i}}))_{n} \is \coprod_{p=0}^{i-1}(\tp{X_{n+p}}{(Y_{\leq 0})_{-p}})=\prod_{p=0}^{i-1}(\tp{X_{n+p}}{(Y_{\leq 0})_{-p}}).$$
This implies that
$${\lim}_{i\in\mathbb{N}}((\tp{X}{Y_{\leq 0}})/(\tp{X}{Y_{\leq -i}}))_{n} \is \prod_{p\in\mathbb{Z}}(\tp{X_{n+p}}{(Y_{\leq 0})_{-p}})=(\btp{X}{Y_{\leq 0}})_{n}.$$
Now it is straightforward to verify
$$\btp{X}{Y_{\leq 0}} \is {\lim}_{i\in\mathbb{N}}((\tp{X}{Y_{\leq 0}})/(\tp{X}{Y_{\leq -i}})).$$

Since $\lim_{i\in\mathbb{N}}(\tp{X}{Y_{\leq -i}})=0$, there is an exact sequence
$$0 \to \tp{X}{Y_{\leq 0}} \to \btp{X}{Y_{\leq 0}} \to {\lim}^{1}_{i\in\mathbb{N}}(\tp{X}{Y_{\leq -i}}) \to 0$$
by Lemma \ref{exact sequence} and the isomorphism proved above. Thus one gets $$\ttp{X}{Y_{\leq 0}} \is {\lim}^{1}_{i\in\mathbb{N}}(\tp{X}{Y_{\leq -i}}).$$

In the general case, where $\sup Y=k\in\ZZ$, we notice that $Y=\susp^{k}(\susp^{-k}Y)_{\leq 0}$ and $(\susp^{-k}Y)_{\leq -i}=\susp^{-k}Y_{\leq k-i}$. Thus one has,
  \begin{align*}
    \btp{X}{Y}
    &=\susp^k(\btp{X}{(\susp^{-k}Y)_{\leq 0}})\\
    &\is \susp^k{\lim}_{i\in\mathbb{N}}((\tp{X}{(\susp^{-k}Y)_{\leq 0}})/(\tp{X}{(\susp^{-k}Y)_{\leq -i}}))\\
    &\is \susp^k{\lim}_{i\in\mathbb{N}}((\tp{X}{\susp^{-k}Y_{\leq k}})/(\tp{X}{\susp^{-k}Y_{\leq k-i}})\\
    &\is {\lim}_{i\in\mathbb{N}}((\tp{X}{Y})/(\tp{X}{Y_{\leq k-i}}))\:,
  \end{align*}
and
  \begin{align*}
    \ttp{X}{Y}
    &=\susp^k(\ttp{X}{(\susp^{-k}Y)_{\leq 0}})\\
    &\is \susp^k{\lim}^{1}_{i\in\mathbb{N}}(\tp{X}{(\susp^{-k}Y)_{\leq -i}})\\
    &\is \susp^k{\lim}^{1}_{i\in\mathbb{N}}(\tp{X}{\susp^{-k}Y_{\leq k-i}})\\
    &\is {\lim}^{1}_{i\in\mathbb{N}}(\tp{X}{Y_{\leq k-i}})\:,
  \end{align*}
as desired.
\end{prf*}

\begin{cor}\label{stensor sequence}
Let $X$ be an $\Rop$-complex and $Y$ a bounded above $R$-complex with $\sup Y=k$. Then there exists an exact sequence
$$0 \to \prod_{i\in\mathbb{N}}(\tp{X}{Y_{\leq k-i}}) \to \prod_{i\in\mathbb{N}}(\tp{X}{Y_{\leq k-i}}) \to \ttp{X}{Y} \to 0.$$
\end{cor}
\begin{prf*}
Since $\lim_{i\in\mathbb{N}}(\tp{X}{Y_{\leq k-i}})=0$ and ${\lim}^{1}_{i\in\mathbb{N}}(\tp{X}{Y_{\leq k-i}}) \is \ttp{X}{Y}$ by Theorem \ref{stensor}, the desired exact sequence now follows from \ref{limits}. We notice that the map from $\prod_{i\in\mathbb{N}}(\tp{X}{Y_{\leq k-i}})$ to $\prod_{i\in\mathbb{N}}(\tp{X}{Y_{\leq k-i}})$ in the statement is $1-\varepsilon$, where $\varepsilon^{uv}: \tp{X}{Y_{\leq k-v}} \to \tp{X}{Y_{\leq k-u}}$ for $u\leq v$ is induced by the filtration
$Y_{\leq k}\supseteq Y_{\leq k-1}\supseteq Y_{\leq k-2}\supseteq \cdots$; see \ref{limits} and \ref{filtration}.
\end{prf*}

\begin{cor}\label{stensor homology}
Let $X$ be an $\Rop$-complex and $Y$ a bounded above $R$-complex with $\sup Y=k$. Then for each $n\in\mathbb{Z}$, there exists an exact sequence
$$0 \to {\lim}^{1}_{i\in\mathbb{N}}\H[n+1]{\tp{X}{Y_{\leq k-i}}} \to \H[n+1]{\ttp{X}{Y}} \to {\lim}_{i\in\mathbb{N}}\H[n]{\tp{X}{Y_{\leq k-i}}} \to 0.$$
In particular, $\H[n+1]{\ttp{X}{Y}}=0$ if and only if ${\lim}^{1}_{i\in\mathbb{N}}\H[n+1]{\tp{X}{Y_{\leq k-i}}}=0={\lim}_{i\in\mathbb{N}}\H[n]{\tp{X}{Y_{\leq k-i}}}$.
\end{cor}
\begin{prf*}
By Corollary \ref{stensor sequence} there is an exact sequence
$$0 \to \prod_{i\in\mathbb{N}}(\tp{X}{Y_{\leq k-i}}) \to \prod_{i\in\mathbb{N}}(\tp{X}{Y_{\leq k-i}}) \to \ttp{X}{Y} \to 0.$$
Thus one gets the following exact sequence
$$\cdots \to \prod_{i\in\mathbb{N}}\H[n+1]{\tp{X}{Y_{\leq k-i}}} \to \prod_{i\in\mathbb{N}}\H[n+1]{\tp{X}{Y_{\leq k-i}}} \to \H[n+1]{\ttp{X}{Y}}$$
$$\to \prod_{i\in\mathbb{N}}\H[n]{\tp{X}{Y_{\leq k-i}}} \to \prod_{i\in\mathbb{N}}\H[n]{\tp{X}{Y_{\leq k-i}}} \to \cdots,$$
which yields the desired exact sequence from the definitions of the lim and lim$^{1}$ groups.
\end{prf*}

\begin{ipg}\label{ML}
Recall that an $\mathbb{N}$-inverse system ${\{\delta_{uv}: M_{v} \to M_{u}}\}_{u\leq v}$ of $R$-modules satisfies the Mittag-Leffler condition if for each $i\in \mathbb{N}$ there exists an index $j\in \mathbb{N}$ with $j\geq i$, such that $\Im\delta_{ij}=\Im\delta_{ik}$ for each $k\in \mathbb{N}$ with $k\geq j$. It is clear that if $\delta_{i,i+1}$ is surjective for each $i\gg 0$ then the $\mathbb{N}$-inverse system ${\{\delta_{uv}: M_{v} \to M_{u}}\}_{u\leq v}$ satisfies the Mittag-Leffler condition. Grothendieck proved in \cite{G} that if the $\mathbb{N}$-inverse system ${\{\delta_{uv}: M_{v} \to M_{u}}\}_{u\leq v}$ satisfies the Mittag-Leffler condition then one has $\lim^1_{i\in\mathbb{N}}M_{i}=0$. Moreover, following \cite[Corollary 6]{E}, $\lim^1_{i\in\mathbb{N}}M_{i}^{(\mathbb{N})}=0$ if and only if the $\mathbb{N}$-inverse system ${\{\delta_{uv}: M_{v} \to M_{u}}\}_{u\leq v}$ satisfies the Mittag-Leffler condition.
\end{ipg}

\begin{cor}\label{stensor ML}
Let $X$ be an $\Rop$-complex and $Y$ a bounded above $R$-complex with $\sup Y=k$, and let $n\in\mathbb{Z}$. If $\H[n]{\ttp{X^{(\mathbb{N})}}{Y}}=0$, then the $\mathbb{N}$-inverse system ${\{\delta_{uv}: \H[n]{\tp{X}{Y_{\leq k-v}}} \to \H[n]{\tp{X}{Y_{\leq k-u}}}}\}_{u\leq v}$ satisfies the Mittag-Leffler condition.
\end{cor}
\begin{prf*}
If $\H[n]{\ttp{X^{(\mathbb{N})}}{Y}}=0$, then by Corollary \ref{stensor homology}, ${\lim}^{1}_{i\in\mathbb{N}}\H[n]{\tp{X^{(\mathbb{N})}}{Y_{\leq k-i}}}=0$, and so one gets ${\lim}^{1}_{i\in\mathbb{N}}(\H[n]{\tp{X}{Y_{\leq k-i}}})^{(\mathbb{N})}=0$, which implies that the $\mathbb{N}$-inverse system ${\{\delta_{uv}: \H[n]{\tp{X}{Y_{\leq k-v}}} \to \H[n]{\tp{X}{Y_{\leq k-u}}}}\}_{u\leq v}$ satisfies the Mittag-Leffler condition; see \ref{ML}.
\end{prf*}

Checking the proof of \cite[Lemma 4.1]{EM}, one gets the following result.

\begin{lem}\label{duality}
Let ${\{\delta_{uv}: X_{v} \to X_{u}}\}_{u\leq v}$ be an $\mathbb{N}$-inverse system of $R$-modules satisfying the Mittag-Leffler condition. If \  $\lim_{i\in\mathbb{N}}X_{i}=0$, then one has
$${\colim}_{i\in\mathbb{N}}\Hom[\mathbb{Z}]{X_{i}}{\mathbb{Q}/\mathbb{Z}}=0.$$
\end{lem}

The next proposition will be used to prove our main results advertised in the introduction.

\begin{prp}\label{colim}
Let $X$ be an $\Rop$-complex and $Y$ a bounded above $R$-complex with $\sup Y=k$, and let $n\in\mathbb{Z}$. If $\H[n]{\ttp{X^{(\mathbb{N})}}{Y}}=0=\H[n+1]{\ttp{X}{Y}}$, then one has $$\colim_{i\in\mathbb{N}}\H[-n]{\Hom[\Rop]{X}{\Hom[\mathbb{Z}]{Y}{\mathbb{Q}/\mathbb{Z}}_{\geq i-k}}}=0$$
 and
 $$\colim_{i\in\mathbb{N}}\H[-n]{\Hom{Y_{\leq k-i}}{\Hom[\mathbb{Z}]{X}{\mathbb{Q}/\mathbb{Z}}}}=0.$$
\end{prp}
\begin{prf*}
The $\mathbb{N}$-inverse system ${\{\delta_{uv}: \H[n]{\tp{X}{Y_{\leq k-v}}} \to \H[n]{\tp{X}{Y_{\leq k-u}}}}\}_{u\leq v}$ satisfies the Mittag-Leffler condition by Corollary \ref{stensor ML}. The vanishing of $\H[n+1]{\ttp{X}{Y}}$ implies that $\lim_{i\in\mathbb{N}}\H[n]{\tp{X}{Y_{\leq k-i}}}=0$; see Corollary \ref{stensor homology}. Thus by Lemma \ref{duality}, one has
 \begin{align*}
    \colim_{i\in\mathbb{N}}\H[-n]{\Hom[\mathbb{Z}]{\tp{X}{Y_{\leq k-i}}}{\mathbb{Q}/\mathbb{Z}}}
    & \is \colim_{i\in\mathbb{N}}\Hom[\mathbb{Z}]{\H[n]{\tp{X}{Y_{\leq k-i}}}}{\mathbb{Q}/\mathbb{Z}}=0.
  \end{align*}
Now the desired equations hold by the adjoint isomorphism.
\end{prf*}

We end this section with the following result that will be used in the next section.

\begin{prp}\label{associativity}
Let $X$ be an $\Rop$-complex, let $Y$ be a bounded $(R,\Sop)$-complex, and let $Z$ be an $S$-complex. Then there is an isomorphism of $\mathbb{Z}$-complexes,
$$\ttp[S]{(\tp{X}{Y})}{Z} \to \ttp{X}{(\tp[S]{Y}{Z})},$$
which is functorial in $X$, $Y$ and $Z$.
\end{prp}
\begin{prf*}
Consider the commutative diagram of $\mathbb{Z}$-complexes:
$$\xymatrix@C=0.4cm{
  0 \ar[r] & \tp[S]{(\tp{X}{Y})}{Z} \ar[d] \ar[r] & \btp[S]{(\tp{X}{Y})}{Z} \ar[d]_{\alpha} \ar[r] & \ttp[S]{(\tp{X}{Y})}{Z} \ar[r] & 0 \\
  0 \ar[r] & \tp{X}{(\tp[S]{Y}{Z})} \ar[r] & \btp{X}{(\tp[S]{Y}{Z})} \ar[r] & \ttp{X}{(\tp[S]{Y}{Z})} \ar[r] & 0.}$$
We notice that $\btp{X}{Y}=\tp{X}{Y}$ and $\btp[S]{Y}{Z}=\tp[S]{Y}{Z}$ since $Y$ is bounded. Then the second vertical map $\alpha$ is an isomorphism by \cite[Proposition A4]{CCLP01}. The first one is clearly an isomorphism. So one gets an isomorphism $$\ttp[S]{(\tp{X}{Y})}{Z} \to \ttp{X}{(\tp[S]{Y}{Z})},$$
which is clearly functorial in $X$, $Y$ and $Z$.
\end{prf*}

\section{stable homology with respect to a semidualizing module} 

\noindent
{\bf Convention.} In this section, $R$ is a commutative noetherian ring, and $C$ is a semidualizing $R$-module.

\begin{dfn}\label{stable homology}
Let $\mathcal{X}$ (resp., $\mathcal{Y}$) be a precovering (resp., preenveloping) class of $R$-modules. For $R$-modules $M$ and $N$, let $X \to M$ be a proper $\mathcal{X}$-resolution of $M$, and $N \to Y$ be a proper $\mathcal{Y}$-coresolution of $N$. For each $n\in\ZZ$, the $n$th \emph{stable homology} of $M$ and $N$ with respect to $\mathcal{X}$ and $\mathcal{Y}$ is
  \begin{equation*}
    \Stor[\mathcal{X}\mathcal{Y}]{n}{M}{N} \deq \H[n+1]{\ttp{X}{Y}}\:.
  \end{equation*}
\end{dfn}

\begin{ipg}
Following \cite[Section 8.2]{EJ}, any two proper $\mathcal{X}$-resolutions of $M$, and similarly any two proper $\mathcal{Y}$-coresolution of $N$, are homotopy equivalent. Thus by \cite[1.5(d)]{CCLP01}, the above definition is independent of the choices of (co)resolutions. We notice that $\Stor[\mathcal{P}\mathcal{I}]{n}{M}{N}$ is the classical stable homology, $\Stor{n}{M}{N}$, of $M$ and $N$ defined by Goichot \cite{Go}; see also \cite{CCLP01}.
\end{ipg}

We denote by $\mathcal{P}_{C}$ (resp., $\mathcal{F}_{C}$, $\mathcal{I}_{C}$) the class of $R$-modules $\tp{C}{P}$ (resp., $\tp{C}{F}$, $\Hom{C}{I}$) with $P$ projective (resp., $F$ flat, $I$ injective). Then $\mathcal{P}_{C}$ and $\mathcal{F}_{C}$ are precovering and $\mathcal{I}_{C}$ is preenveloping; see e.g. Holm and White \cite[Proposition 5.3]{HW}. In the next lemma, (a) and (b) can be found in \cite[Lemma 3.1]{SSWTY}, (c) can be proved as in \cite[Lemma 3.1(c)]{SSWTY}, and (d) is from \cite[Lemma 2.1(b)]{TW}.

\begin{lem}\label{proper resolution}
Let $M$ be an $R$-module.
  \begin{eqc}
  \item[$\mathrm{(a)}$] If $F\qra\Hom{C}{M}$ is a proper flat (resp., projective) resolution, then $\tp{C}{F}\to M$ is a proper $\mathcal{F}_{C}$ (resp., $\mathcal{P}_{C}$)-resolution of $M$.
  \item[$\mathrm{(b)}$] If $G \to M$ is a proper $\mathcal{F}_{C}$ (resp., $\mathcal{P}_{C}$)-resolution of $M$, then $\Hom{C}{G}\qra\Hom{C}{M}$ is a proper flat (resp., projective)-resolution of $\Hom{C}{M}$.
  \item[$\mathrm{(c)}$] If $\tp{C}{M} \qra I$ is an injective resolution of $\tp{C}{M}$, then $M \to \Hom{C}{I}$ is a proper $\mathcal{I}_C$-coresolution.
  \item[$\mathrm{(d)}$] If $M\to J$ is a proper $\mathcal{I}_{C}$-coresolution of $M$, then $\tp{C}{M}\qra\tp{C}{J}$ is an injective resolution of $\tp{C}{M}$.
  \end{eqc}
\end{lem}

\begin{prp}\label{Cflat resolution}
Let $M$ and $N$ be $R$-modules. Then there are isomorphisms
$$\Stor[\mathcal{P}_C\mathcal{I}_C]{n}{M}{N} \is \Stor{n}{\Hom{C}{M}}{\tp{C}{N}} \is \Stor[\mathcal{F}_C\mathcal{I}_C]{n}{M}{N},$$
which are functorial in $M$ and $N$.
\end{prp}
\begin{prf*}
Let $P \qra \Hom{C}{M}$ be a projective resolution of $\Hom{C}{M}$, and let $\tp{C}{N} \qra I$ be an injective resolution of $\tp{C}{N}$. Then by Lemma \ref{proper resolution}(a)(c), $\tp{C}{P} \to M$ is a proper $\mathcal{P}_C$-resolution of $M$, and $N \to \Hom{C}{I}$ is a proper $\mathcal{I}_C$-coresolution, and so one gets
 \begin{align*}
    \Stor[\mathcal{P}_C\mathcal{I}_C]{n}{M}{N}
    & = \H[n+1]{\ttp{(\tp{C}{P})}{\Hom{C}{I}}}\\
    & \is \H[n+1]{\ttp{P}{(\tp{C}{\Hom{C}{I}})}}\\
    & \is \H[n+1]{\ttp{P}{I}}\\
    & \is \Stor{n}{\Hom{C}{M}}{\tp{C}{N}},
  \end{align*}
where the first isomorphism follows from Proposition \ref{associativity}, and the second one holds since $I$ is a complex of injective $R$-modules.

The isomorphism $\Stor[\mathcal{F}_C\mathcal{I}_C]{n}{M}{N} \is \Stor{n}{\Hom{C}{M}}{\tp{C}{N}}$ can be proved similarly by taking a proper flat resolution $F \qra \Hom{C}{M}$ and using Lemma \ref{proper resolution}(a) and \cite[Proposition 2.6]{CCLP01}.

Now it is straightforward to verify that the desired isomorphisms are functorial in $M$ and $N$.
\end{prf*}

\begin{lem}\label{dimension shift}
Let $M$ be an $R$-module and let $n\in\ZZ$.
  \begin{eqc}
  \item[$\mathrm{(a)}$] If $\Stor[\mathcal{P}_C\mathcal{I}_C]{n-1}{-}{M}=0$, then $\Stor[\mathcal{P}_C\mathcal{I}_C]{n}{-}{M}=0$.
  \item[$\mathrm{(b)}$] If $\Stor[\mathcal{P}_C\mathcal{I}_C]{n+1}{M}{-}=0$, then $\Stor[\mathcal{P}_C\mathcal{I}_C]{n}{M}{-}=0$.
  \end{eqc}
\end{lem}
\begin{prf*}
(a) For an $R$-module $M'$, by \cite[Proposition 5.3(b)]{HW} there is a complex $0 \to K \to P \to M' \to 0$ with $P\in\mathcal{P}_{C}$ such that the sequence
$$0 \to \Hom{P'}{K} \to \Hom{P'}{P} \to \Hom{P'}{M'} \to 0$$
is exact for each $P'\in\mathcal{P}_{C}$. In particular, the sequence
$$0 \to \Hom{C}{K} \to \Hom{C}{P} \to \Hom{C}{M'} \to 0$$
is exact. Since $\Hom{C}{P}$ is projective, one gets
$$\Stor{n}{\Hom{C}{M'}}{\tp{C}{M}} \is \Stor{n-1}{\Hom{C}{K}}{\tp{C}{M}},$$
and so by Proposition \ref{Cflat resolution}, $\Stor[\mathcal{P}_C\mathcal{I}_C]{n}{M'}{M} \is \Stor[\mathcal{P}_C\mathcal{I}_C]{n-1}{K}{M}=0$, which yields $\Stor[\mathcal{P}_C\mathcal{I}_C]{n}{-}{M}=0$.

(b) Let $N$ be an $R$-module. Then by \cite[Proposition 5.3(c)]{HW} there is a complex $0 \to N \to I \to K \to 0$ with $I\in\mathcal{I}_{C}$ such that the sequence
$$0 \to \Hom{K}{I'} \to \Hom{I}{I'} \to \Hom{N}{I'} \to 0$$
is exact for each $I'\in\mathcal{I}_{C}$. Since $C^{\vee}=\Hom[\mathbb{Z}]{C}{\mathbb{Q}/\mathbb{Z}}$ is in $\mathcal{I}_{C}$, the sequence
$$0 \to \Hom{K}{C^{\vee}} \to \Hom{I}{C^{\vee}} \to \Hom{N}{C^{\vee}} \to 0$$
is exact, which implies that the sequence
$$0 \to \tp{C}{N} \to \tp{C}{I} \to \tp{C}{K} \to 0$$
is exact. We notice that $\tp{C}{I}$ is injective. Thus one gets
$$\Stor{n}{\Hom{C}{M}}{\tp{C}{N}} \is \Stor{n+1}{\Hom{C}{M}}{\tp{C}{K}},$$
and so by Proposition \ref{Cflat resolution}, $\Stor[\mathcal{P}_C\mathcal{I}_C]{n}{M}{N} \is \Stor[\mathcal{P}_C\mathcal{I}_C]{n+1}{M}{K}=0$, which yields $\Stor[\mathcal{P}_C\mathcal{I}_C]{n}{M}{-}=0$.
\end{prf*}

Now we are in a position to give the main results of this section described in the introduction.

\begin{thm}\label{cinj dim-Tor}
For an $R$-module $N$, the following conditions are equivalent:
  \begin{eqc}
  \item $\mathcal{I}_C{\text-}\id_{R}N<\infty$.
  \item $\Stor[\mathcal{P}_C\mathcal{I}_C]{n}{-}{N}=0$ for each $n\in\ZZ$.
  \item $\Stor[\mathcal{P}_C\mathcal{I}_C]{n}{-}{N}=0$ for some $n<0$.
  \end{eqc}
\end{thm}
\begin{prf*}
\proofofimp{i}{ii} Since $\mathcal{I}_C{\text-}\id_{R}N<\infty$, there is a proper $\mathcal{I}_{C}$-coresolution $N \to I$ with $I$ bounded. Thus for each $R$-module $M$ with $P \to M$ a proper $\mathcal{P}_{C}$-resolution, one has $\Stor[\mathcal{P}_C\mathcal{I}_C]{n}{M}{N}=\H[n+1]{\ttp{P}{I}}=0$.

\eqclbl{ii}$\implies$\eqclbl{iii} is clear.

\proofofimp{iii}{i} We first notice that $\Stor[\mathcal{P}_C\mathcal{I}_C]{0}{-}{N}=0=\Stor[\mathcal{P}_C\mathcal{I}_C]{-1}{-}{N}=0$ by Lemma \ref{dimension shift}.

Let $M$ be an $R$-module and $F \qra \Hom{C}{M}$ a proper flat resolution. Then by Lemma \ref{proper resolution}(a), $\tp{C}{F} \to M$ is a proper $\mathcal{F}_C$-resolution of $M$. Let $N \to I$ be a proper $\mathcal{I}_C$-coresolution of $N$. Since $\Stor[\mathcal{F}_C\mathcal{I}_C]{0}{M}{N} \is \Stor[\mathcal{P}_C\mathcal{I}_C]{0}{M}{N}=0$ by Proposition \ref{Cflat resolution}, one gets $\H[1]{\ttp{(\tp{C}{F})}{I}}=0$.

On the other hand, one has $\Stor[\mathcal{P}_C\mathcal{I}_C]{-1}{M^{(\mathbb{N})}}{N}=0$, so by Proposition \ref{Cflat resolution},
$$\Stor{-1}{(\Hom{C}{M})^{(\mathbb{N})}}{\tp{C}{N}}=0.$$
Note that $F \qra \Hom{C}{M}$ is a flat resolution, and so $F^{(\mathbb{N})} \qra (\Hom{C}{M})^{(\mathbb{N})}$ is a flat resolution of $(\Hom{C}{M})^{(\mathbb{N})}$. Since $\tp{C}{N} \qra \tp{C}{I}$ is an injective resolution by Lemma \ref{proper resolution}(d), one gets $\H[0]{\ttp{F^{(\mathbb{N})}}{(\tp{C}{I})}}=0$; see \cite[Proposition 2.6]{CCLP01}. Thus we have
 \begin{align*}
    \H[0]{\ttp{(\tp{C}{F})^{(\mathbb{N})}}{I}}
    & \is \H[0]{\ttp{(\tp{C}{F^{(\mathbb{N})})}}{I}} \is \H[0]{\ttp{F^{(\mathbb{N})}}{(\tp{C}{I})}}=0,
  \end{align*}
where the second isomorphism follows from Proposition \ref{associativity}.

Now by Proposition \ref{colim} one gets
$$\colim_{i\in\mathbb{N}}\H[0]{\Hom{\tp{C}{F}}{\Hom[\mathbb{Z}]{I}{\mathbb{Q}/\mathbb{Z}}_{\geq i}}}=0.$$
We notice that $\tp{C}{F} \to M$ is a proper $\mathcal{F}_C$-resolution of $M$, and $\Hom[\mathbb{Z}]{I}{\mathbb{Q}/\mathbb{Z}} \to \Hom[\mathbb{Z}]{N}{\mathbb{Q}/\mathbb{Z}}$ is a proper $\mathcal{F}_C$-resolution of $\Hom[\mathbb{Z}]{N}{\mathbb{Q}/\mathbb{Z}}$. Then by Proposition \ref{csExt colim} one gets
$$\sExt[\mathcal{F}_C]{0}{M}{\Hom[\mathbb{Z}]{N}{\mathbb{Q}/\mathbb{Z}}}=0$$
for each $R$-module $M$. Thus $\mathcal{F}_C$-$\pd_{R}\Hom[\mathbb{Z}]{N}{\mathbb{Q}/\mathbb{Z}}<\infty$ by Proposition \ref{Cflat dimension}, and so $\mathcal{I}_C{\text-}\id_{R}{N}<\infty$; see Sather-Wagstaff, Sharif and White \cite[Lemma 4.2]{SWSW}.
\end{prf*}

\begin{thm}\label{cflat dim-Tor}
For an $R$-module $M$, the following conditions are equivalent:
  \begin{eqc}
  \item $\mathcal{F}_C{\text-}\pd_{R}M<\infty$.
  \item $\Stor[\mathcal{P}_C\mathcal{I}_C]{n}{M}{-}=0$ for each $n\in\ZZ$.
  \item $\Stor[\mathcal{P}_C\mathcal{I}_C]{n}{M}{-}=0$ for some $n\geq0$.
  \end{eqc}
Moreover, if $M$ is finitely generated, then \eqclbl{i}--\eqclbl{iii} are equivalent to
  \begin{eqc}
  \item[\eqclbl{i'}] $\mathcal{P}_C{\text-}\pd_{R}M<\infty$.
  \end{eqc}
\end{thm}
\begin{prf*}
\proofofimp{i}{ii} Since $\mathcal{F}_C{\text-}\pd_{R}M<\infty$, there is a proper $\mathcal{F}_{C}$-resolution $F \to M$ with $F$ bounded. Thus for each $R$-module $N$ with $N \to I$ a proper $\mathcal{I}_{C}$-coresolution, one has $\Stor[\mathcal{P}_C\mathcal{I}_C]{n}{M}{N} \is \Stor[\mathcal{F}_C\mathcal{I}_C]{n}{M}{N}=\H[n+1]{\ttp{F}{I}}=0$ by Proposition \ref{Cflat resolution}.

\eqclbl{ii}$\implies$\eqclbl{iii} is clear.

\proofofimp{iii}{i} We first notice that $\Stor[\mathcal{P}_C\mathcal{I}_C]{0}{M}{-}=0=\Stor[\mathcal{P}_C\mathcal{I}_C]{-1}{M}{-}=0$ by Lemma \ref{dimension shift}.

Let $F \qra \Hom{C}{M}$ be a proper flat resolution of $\Hom{C}{M}$. Then
$\tp{C}{F} \to M$ is a proper $\mathcal{F}_C$-resolution by Lemma \ref{proper resolution}(a), and
$$\tp{C}{\Hom[\mathbb{Z}]{M}{\mathbb{Q}/\mathbb{Z}}} \is \Hom[\mathbb{Z}]{\Hom{C}{M}}{\mathbb{Q}/\mathbb{Z}} \qra \Hom[\mathbb{Z}]{F}{\mathbb{Q}/\mathbb{Z}}$$
is an injective resolution of $\tp{C}{\Hom[\mathbb{Z}]{M}{\mathbb{Q}/\mathbb{Z}}}$, and so $$\Hom[\mathbb{Z}]{M}{\mathbb{Q}/\mathbb{Z}} \to \Hom{C}{\Hom[\mathbb{Z}]{F}{\mathbb{Q}/\mathbb{Z}}} \is \Hom[\mathbb{Z}]{\tp{C}{F}}{\mathbb{Q}/\mathbb{Z}}$$
is a proper $\mathcal{I}_C$-coresolution of $\Hom[\mathbb{Z}]{M}{\mathbb{Q}/\mathbb{Z}}$ by Lemma \ref{proper resolution}(c).

Let $N$ be $R$-module, and let $\tp{C}{N} \qra I$ be an injective resolution of $\tp{C}{N}$. Then $N \to \Hom{C}{I}$ is a proper $\mathcal{I}_C$-coresolution by Lemma \ref{proper resolution}(c), and
$$\tp{C}{N^{(\mathbb{N})}} \is (\tp{C}{N})^{(\mathbb{N})} \qra I^{(\mathbb{N})}$$
is an injective resolution of $\tp{C}{N^{(\mathbb{N})}}$, and so
$$N^{(\mathbb{N})} \to \Hom{C}{I^{(\mathbb{N})}} \is (\Hom{C}{I})^{(\mathbb{N})}$$
is a proper $\mathcal{I}_C$-coresolution by Lemma \ref{proper resolution}(c).

Since $\Stor[\mathcal{F}_C\mathcal{I}_C]{0}{M}{N}=0=
\Stor[\mathcal{F}_C\mathcal{I}_C]{-1}{M}{N^{(\mathbb{N})}}$ by Proposition \ref{Cflat resolution}, one gets
$\H[1]{\ttp{(\tp{C}{F})}{\Hom{C}{I}}}=0$, and $$\H[0]{\ttp{(\tp{C}{F})^{(\mathbb{N})}}{\Hom{C}{I}}} \is \H[0]{\ttp{(\tp{C}{F})}{(\Hom{C}{I})^{(\mathbb{N})}}}=0$$
by Proposition \ref{associativity}. Now using Proposition \ref{colim}, one gets
$$\colim_{i\in\mathbb{N}}\H[0]{\Hom{\Hom{C}{I}_{\leq -i}}{\Hom[\mathbb{Z}]{\tp{C}{F}}{\mathbb{Q}/\mathbb{Z}}}}=0,$$
and so $\sExt[\mathcal{I}_C]{0}{N}{\Hom[\mathbb{Z}]{M}{\mathbb{Q}/\mathbb{Z}}}=0$ for each $R$-module $N$ by Proposition \ref{csExt colim-dual}. Thus $\mathcal{I}_C$-$\id_{R}\Hom[\mathbb{Z}]{M}{\mathbb{Q}/\mathbb{Z}}<\infty$ by Proposition \ref{Cinj dimension}, and so $\mathcal{F}_C{\text-}\pd_{R}M<\infty$; see \cite[Lemma 4.2]{SWSW}.

Finally, if $M$ is finitely generated, then by \cite[Theorem 5.5]{SSWTY} the conditions $(i)$ and $(i')$ are equivalent.
\end{prf*}

As a corollary of the above theorems, we give a balance result for stable homology with respect to a semidualizing module.

\begin{cor}\label{c trival}
The following conditions are equivalent for a local ring $R$:
  \begin{eqc}
  \item $\Stor[\mathcal{P}_C\mathcal{I}_C]{n}{M}{N} \is \Stor[\mathcal{P}_C\mathcal{I}_C]{n}{N}{M}$ for all $R$-modules $M$ and $N$, and for each $n\in\ZZ$.
  \item $\mathcal{I}_C$-$\id_{R}C<\infty$.
  \item $C \is R$ and $R$ is Gorenstein.
  \end{eqc}
\end{cor}
\begin{prf*}
\proofofimp{i}{ii} Since $C$ is $C$-projective, one gets $$\Stor[\mathcal{P}_C\mathcal{I}_C]{n}{M}{C} \is \Stor[\mathcal{P}_C\mathcal{I}_C]{n}{C}{M}=0$$
for all $R$-modules $M$ and for each $n\in\ZZ$, and so $\mathcal{I}_C{\text-}\id(C)<\infty$ by Theorem \ref{cinj dim-Tor}.

The implication \eqclbl{ii}$\implies$\eqclbl{iii} follows from Sather-Wagstaff and Yassemi \cite[Lemma 2.11]{SWY}, and \eqclbl{iii}$\implies$\eqclbl{i} holds by \cite[Corollary 4.7]{CCLP01}.
\end{prf*}

\appendix
\section*{Appendix. Stable cohomology}
\stepcounter{section}

\noindent
The next definitions of bounded and stable Hom-complexes can be found in Avramov and Veliche \cite{AV}, and \cite{Go}.

\begin{ipg}
For $R$-complexes $X$ and $Y$, the \emph{bounded Hom-complex} $\bHom{X}{Y}$ is the subcomplex of $\Hom{X}{Y}$ with degree-$n$ term $$\bHom{X}{Y}_{n}=\coprod_{i\in\mathbb{Z}}\Hom{X_{i}}{Y_{n+i}}.$$
We denote by $\SHom{X}{Y}$ the quotient complex $\Hom{X}{Y}/\bHom{X}{Y}$, which is called the \emph{stable Hom-complex}.
\end{ipg}

\begin{prp}\label{adjointness}
Let $X$ and $Z$ be an $R$-complex and an $S$-complex, respectively, and let $Y$ be a bounded ($S, \Rop$)-complex. Then there are isomorphisms of $\mathbb{Z}$-complexes,
$$\bHom[S]{\tp{Y}{X}}{Z} \is \bHom[R]{X}{\Hom[S]{Y}{Z}}$$
and
$$\SHom[S]{\tp{Y}{X}}{Z} \is \SHom[R]{X}{\Hom[S]{Y}{Z}},$$
which are functorial in $X, Y$ and $Z$.
\end{prp}
\begin{prf*}
  For every $n\in\ZZ$ one has,
  \begin{align*}
    \bHom[S]{\tp{Y}{X}}{Z}_{n}
    &=\coprod_{h\in\ZZ}\Hom[S]{(\tp{Y}{X})_{h}}{Z_{n+h}}\\
    &=\coprod_{h\in\ZZ}\Hom[S]{\coprod_{q\in\ZZ}(\tp{Y_{q}}{X_{h-q}})}{Z_{n+h}}\\
    &\is\coprod_{h\in\ZZ}\coprod_{q\in\ZZ}\Hom[S]{\tp{Y_{q}}{X_{h-q}}}{Z_{n+h}}\\
    &=\coprod_{p\in\ZZ}\coprod_{q\in\ZZ}\Hom[S]{\tp{Y_{q}}{X_{p}}}{Z_{n+p+q}}\:.
  \end{align*}
On the other hand, for every $n\in\ZZ$ one has,
  \begin{align*}
    \bHom[R]{X}{\Hom[S]{Y}{Z}}_{n}
    &=\coprod_{p\in\ZZ}\Hom{X_{p}}{\Hom[S]{Y}{Z}_{n+p}}\\
    &=\coprod_{p\in\ZZ}\Hom{X_{p}}{\prod_{q\in\ZZ}\Hom[S]{Y_{q}}{Z_{n+p+q}}}\\
    &\is\coprod_{p\in\ZZ}\coprod_{q\in\ZZ}\Hom{X_{p}}{\Hom[S]{Y_{q}}{Z_{n+p+q}}}\:.
  \end{align*}
Here the isomorphisms in the above computations hold since $Y$ is bounded.

We notice that there is a natural isomorphism of $\ZZ$-modules
\begin{equation*}
    \dmapdef{\rho_{Y_{q}X_{p}Z_{n+p+q}}}{\Hom[S]{\tp{Y_{q}}{X_{p}}}{Z_{n+p+q}}}
    {\Hom{X_{p}}{\Hom[S]{Y_{q}}{Z_{n+p+q}}}}.
  \end{equation*}
Thus one gets an isomorphism of $\mathbb{Z}$-complexes
\begin{equation*}
    \dmapdef{\rho_{YXZ}}{\bHom[S]{\tp{Y}{X}}{Z}}
    {\bHom[R]{X}{\Hom[S]{Y}{Z}}}.
  \end{equation*}
  It is straightforward to verify that $\rho_{YXZ}$ is functorial in $X$, $Y$ and $Z$.

For the second isomorphism in the statement, consider the following commutative diagram of $\mathbb{Z}$-complexes:
$$\xymatrix@C=0.15cm{
  0 \ar[r] & \bHom[S]{\tp{Y}{X}}{Z} \ar[d]_{\rho} \ar[r] & \Hom[S]{\tp{Y}{X}}{Z} \ar[d]_{\varrho} \ar[r] & \SHom[S]{\tp{Y}{X}}{Z} \ar[r] & 0 \\
  0 \ar[r] & \bHom[R]{X}{\Hom[S]{Y}{Z}} \ar[r] & \Hom[R]{X}{\Hom[S]{Y}{Z}} \ar[r] & \SHom[R]{X}{\Hom[S]{Y}{Z}} \ar[r] & 0.}$$
Since $\rho$ and $\varrho$ are isomorphisms, one gets an isomorphism
$$\SHom[S]{\tp{Y}{X}}{Z} \to \SHom[R]{X}{\Hom[S]{Y}{Z}},$$
which is clearly functorial in $X$, $Y$ and $Z$.
\end{prf*}

\begin{ipg}
Let $\mathcal{X}$ be a precovering class of $R$-modules, and let $X_{M} \to M$ and $X_{N} \to N$ be proper $\mathcal{X}$-resolutions of $R$-modules $M$ and $N$, respectively. For each $n\in\ZZ$, the $n$th \emph{stable cohomology} of $M$ and $N$ with respect to $\mathcal{X}$ is
  \begin{equation*}
    \sExt[\mathcal{X}]{n}{M}{N} \deq \H[-n]{\sHom{X_{M}}{X_{N}}}\:.
  \end{equation*}

Dually, let $\mathcal{Y}$ be a preenveloping class of $R$-modules, and let $M \to Y_{M}$ and $N \to Y_{N}$ be proper $\mathcal{Y}$-coresolutions of $M$ and $N$, respectively. For each $n\in\ZZ$, the $n$th \emph{stable cohomology} of $M$ and $N$ with respect to $\mathcal{Y}$ is
  \begin{equation*}
    \sExt[\mathcal{Y}]{n}{M}{N} \deq \H[-n]{\sHom{Y_{M}}{Y_{N}}}\:.
  \end{equation*}
\end{ipg}

\begin{ipg}
Any two proper $\mathcal{X}$-resolutions of $M$, and similarly any two proper $\mathcal{Y}$-coresolutions of $N$, are homotopy equivalent; see \cite[Section 8.2]{EJ}. Thus the above definitions are independent of the choices of (co)resolutions. We notice that $\sExt[\mathcal{P}]{n}{M}{N}$ is the classical stable cohomology, $\sExt{n}{M}{N}$, of $M$ and $N$; see \cite{AV} and \cite{Go}. Also $\sExt[\mathcal{I}]{n}{M}{N}$ is the cohomology given by Nucinkis \cite{N}.
\end{ipg}

\subsection*{\em Stable cohomology with respect to proper flat (injective) resolutions} The proof of the next result can be modelled along the argument in the proof of \cite[Proposition 2.2]{AV}, when the argument is applied to the functor $\Ext[\mathcal{F}]{i}{M}{-}$ that is computed by $\H[-i]{\Hom{F}{-}}$, where $F\qra M$ is a proper flat resolution.
\begin{prp}\label{flat dimension}
For an $R$-module $M$, the following conditions are equivalent:
  \begin{eqc}
  \item $\fd_{R}{M}<\infty$.
  \item $\sExt[\mathcal{F}]{n}{M}{-}=0=\sExt[\mathcal{F}]{n}{-}{M}$ for each $n\in\mathbb{Z}$.
  \item $\sExt[\mathcal{F}]{0}{M}{M}=0$.
  \end{eqc}
\end{prp}

Dually, we have the next result that is proved by Nucinkis in \cite[Theorem 3.7]{N}.

\begin{prp}\label{inj dimension}
For an $R$-module $N$, the following conditions are equivalent:
  \begin{eqc}
  \item $\id_{R}{N}<\infty$.
  \item $\sExt[\mathcal{I}]{n}{N}{-}=0=\sExt[\mathcal{I}]{n}{-}{N}$ for each $n\in\mathbb{Z}$.
  \item $\sExt[\mathcal{I}]{0}{N}{N}=0$.
  \end{eqc}
\end{prp}

\begin{prp}\label{sExt colim}
Let $M$ and $N$ be $R$-modules with proper flat resolutions $F \qra
  M$ and $F' \qra N$, respectively. For every $n\in\ZZ$ there is an isomorphism
  \begin{equation*}
    \sExt[\mathcal{F}]{n}{M}{N} \dis \colim_{i\in\mathbb{N}}\H[-n]{\Hom{F}{F'_{\geq i}}}.
  \end{equation*}
\end{prp}
\begin{prf*}
Set $\Omega_{s}M=\Coker(F_{s+1}\to F_{s})$ and $\Omega_{s}N=\Coker(F'_{s+1}\to F'_{s})$. Using a similar proof as proved in \cite[Theorem 3.6]{N}, one gets a natural isomorphism
  \begin{equation*}
    \colim_{i\in\mathbb{N}}\Ext[\mathcal{F}]{i}{M}{\Omega_{i-n}N} \dis \colim_{i\in\mathbb{N}}\Hom{\Omega_{i}M}{\Omega_{i-n}N}/\FHom{\Omega_{i}M}{\Omega_{i-n}N}.
  \end{equation*}
Here $\FHom{\Omega_{i}M}{\Omega_{i-n}N}$ denotes the set of all homomorphisms of $R$-modules $f\in\Hom{\Omega_{i}M}{\Omega_{i-n}N}$ factoring through a flat $R$-module. As proved in \cite[Theorem 4.4]{N} (see also \cite[B.2]{CCLP02}), one gets an isomorphism
  \begin{equation*}
    \sExt[\mathcal{F}]{n}{M}{N} \dis \colim_{i\in\mathbb{N}}\Hom{\Omega_{i}M}{\Omega_{i-n}N}/\FHom{\Omega_{i}M}{\Omega_{i-n}N}.
  \end{equation*}
On the other hand, we notice that $\susp^{-i}F'_{\geq i} \qra \Omega_{i}N$ is a proper flat resolution. Thus one has,
 \begin{align*}
    \colim_{i\in\mathbb{N}}\Ext[\mathcal{F}]{i}{M}{\Omega_{i-n}N}
    & \is \colim_{i\in\mathbb{N}}\Ext[\mathcal{F}]{i+n}{M}{\Omega_{i}N}\\
    & \is \colim_{i\in\mathbb{N}}\H[-i-n]{\Hom{F}{\susp^{-i}F'_{\geq i}}}\\
    & \is \colim_{i\in\mathbb{N}}\H[-n]{\Hom{F}{F'_{\geq i}}}\:,
 \end{align*}
where the second isomorphism follows from Christensen, Frankild and Holm \cite[Proposition 2.6]{CFH}. Now one gets the isomorphism in the statement.
\end{prf*}

Dually, one gets the following result, which is proved in \cite[Proposition 1.1(iii)]{EM}.

\begin{prp}\label{sExt colim-dual}
Let $M$ and $N$ be $R$-modules with injective resolutions $M \qra
  I$ and $N \qra I'$, respectively. For every $n\in\ZZ$ there is an isomorphism
  \begin{equation*}
    \sExt[\mathcal{I}]{n}{M}{N} \dis \colim_{i\in\mathbb{N}}\H[-n]{\Hom{I_{\leq -i}}{I'}}.
  \end{equation*}
\end{prp}

\subsection*{\em Stable cohomology with respect to a semidualizing module} In this subsection, we assume that $R$ is a commutative noetherian ring, and let $C$ be a semidualizing $R$-module.

\begin{lem}\label{sExt}
Let $M$ and $N$ be $R$-modules. Then there is an isomorphism
$$\sExt[\mathcal{F}_{C}]{n}{M}{N} \is \sExt[\mathcal{F}]{n}{\Hom{C}{M}}{\Hom{C}{N}},$$
which is functorial in $M$ and $N$.
\end{lem}
\begin{prf*}
Let $F \qra \Hom{C}{M}$ and $F' \qra \Hom{C}{N}$ be proper flat resolutions of $\Hom{C}{M}$ and $\Hom{C}{N}$, respectively. Then by Lemma \ref{proper resolution}(a), $\tp{C}{F} \to M$ and $\tp{C}{F'} \to N$ are proper $\mathcal{F}_{C}$-resolutions of $M$ and $N$, respectively. Thus one has,
 \begin{align*}
    \sExt[\mathcal{F}_{C}]{n}{M}{N}
    &= \H[-n]{\SHom{\tp{C}{F}}{\tp{C}{F'}}}\\
    & \is \H[-n]{\SHom{F}{\Hom{C}{\tp{C}{F'}}}}\\
    & \is \H[-n]{\SHom{F}{F'}}\\
    & \is \sExt[\mathcal{F}]{n}{\Hom{C}{M}}{\Hom{C}{N}}\:,
 \end{align*}
where the first isomorphism follows from Proposition \ref{adjointness}, and the second one holds since $F'$ is a complex of flat $R$-modules. It is straightforward to verify that the desired isomorphism is functorial in $M$ and $N$.
\end{prf*}

The next result can be proved dually using Lemma \ref{proper resolution}(c) and Proposition \ref{adjointness}.

\begin{lem}\label{sExt-dual}
Let $M$ and $N$ be $R$-modules. Then there is an isomorphism
$$\sExt[\mathcal{I}_{C}]{n}{M}{N} \is \sExt[\mathcal{I}]{n}{\tp{C}{M}}{\tp{C}{N}},$$
which is functorial in $M$ and $N$.
\end{lem}

\begin{prp}\label{Cflat dimension}
For an $R$-module $M$, the following conditions are equivalent.
  \begin{eqc}
  \item $\mathcal{F}_{C}{\text-}\pd_{R}{M}<\infty$.
  \item $\sExt[\mathcal{F}_{C}]{n}{M}{-}=0=\sExt[\mathcal{F}_{C}]{n}{-}{M}$ for each $n\in\mathbb{Z}$.
  \item $\sExt[\mathcal{F}_{C}]{0}{M}{M}=0$.
  \end{eqc}
\end{prp}
\begin{prf*}
\proofofimp{i}{ii} Since $\mathcal{F}_{C}{\text-}\pd_{R}{M}<\infty$, there is a proper $\mathcal{F}_{C}$-resolution $F \to M$ with $F$ bounded, and so $\SHom{F}{-}=0=\SHom{-}{F}$. Thus one gets $\sExt[\mathcal{F}_{C}]{n}{M}{-}=0=\sExt[\mathcal{F}_{C}]{n}{-}{M}$ for each $n\in\mathbb{Z}$.

\eqclbl{ii}$\implies$\eqclbl{iii} is clear.

\proofofimp{iii}{i} By Proposition \ref{sExt}, one gets $\sExt[\mathcal{F}]{0}{\Hom{C}{M}}{\Hom{C}{M}} \is \sExt[\mathcal{F}_{C}]{0}{M}{M}=0$, and so $\fd_{R}\Hom{C}{M}<\infty$ by Proposition \ref{flat dimension}. Thus one gets $\mathcal{F}_{C}{\text-}\pd_{R}{M}<\infty$; see \cite[Proposition 5.2(b)]{SSWTY}.
\end{prf*}

The next result can be proved dually using Propositions \ref{inj dimension} and \ref{sExt-dual}, and \cite[Theorem 2.11(b)]{TW}.

\begin{prp}\label{Cinj dimension}
For an $R$-module $N$, the following conditions are equivalent.
  \begin{eqc}
  \item $\mathcal{I}_{C}{\text-}\id_{R}{N}<\infty$.
  \item $\sExt[\mathcal{I}_{C}]{n}{N}{-}=0=\sExt[\mathcal{I}_{C}]{n}{-}{N}$ for each $n\in\mathbb{Z}$.
  \item $\sExt[\mathcal{I}_{C}]{0}{N}{N}=0$.
  \end{eqc}
\end{prp}

\begin{prp}\label{csExt colim}
Let $M$ and $N$ be $R$-modules with proper $\mathcal{F}_{C}$-resolutions $F \to
  M$ and $F' \to N$, respectively. For every $n\in\ZZ$ there is an isomorphism
  \begin{equation*}
    \sExt[\mathcal{F}_{C}]{n}{M}{N} \dis \colim_{i\in\mathbb{N}}\H[-n]{\Hom{F}{F'_{\geq i}}}.
  \end{equation*}
\end{prp}
\begin{prf*}
By Lemma \ref{proper resolution}(b), $\Hom{C}{F} \qra \Hom{C}{M}$ and $\Hom{C}{F'} \qra \Hom{C}{N}$ are proper flat resolutions of $\Hom{C}{M}$ and $\Hom{C}{N}$, respectively. Thus we have,
 \begin{align*}
    \sExt[\mathcal{F}_{C}]{n}{M}{N}
    &\is \sExt[\mathcal{F}]{n}{\Hom{C}{M}}{\Hom{C}{N}}\\
    & \is \colim_{i\in\mathbb{N}}\H[-n]{\Hom{\Hom{C}{F}}{\Hom{C}{F'}_{\geq i}}}\\
    & = \colim_{i\in\mathbb{N}}\H[-n]{\Hom{\Hom{C}{F}}{\Hom{C}{F'_{\geq i}}}}\\
    & \is \colim_{i\in\mathbb{N}}\H[-n]{\Hom{\tp{C}{\Hom{C}{F}}}{F'_{\geq i}}}\\
    & \is \colim_{i\in\mathbb{N}}\H[-n]{\Hom{F}{F'_{\geq i}}}\:,
 \end{align*}
where the first isomorphism follows from Proposition \ref{sExt}, the second one follows from Proposition \ref{sExt colim}, and the last one holds since $F$ is a complex of $C$-flat $R$-modules.
\end{prf*}

Dually, we have the following result.

\begin{prp}\label{csExt colim-dual}
Let $M$ and $N$ be $R$-modules with proper $\mathcal{I}_{C}$-coresolutions $M \to
  I$ and $N \to I'$, respectively. For every $n\in\ZZ$ there is an isomorphism
  \begin{equation*}
    \sExt[\mathcal{I}_{C}]{n}{M}{N} \dis \colim_{i\in\mathbb{N}}\H[-n]{\Hom{I_{\leq -i}}{I'}}.
  \end{equation*}
\end{prp}

\section*{Acknowledgments}

\noindent
We thank Ioannis Emmanouil and Panagiota Manousaki for making \cite{EM} available to us and for discussions regarding this work. We also thank the anonymous referees for several corrections and valuable comments that improved the presentation at several points.

\def\cprime{$'$}
  \providecommand{\arxiv}[2][AC]{\mbox{\href{http://arxiv.org/abs/#2}{\sf
  arXiv:#2 [math.#1]}}}
  \providecommand{\oldarxiv}[2][AC]{\mbox{\href{http://arxiv.org/abs/math/#2}{\sf
  arXiv:math/#2
  [math.#1]}}}\providecommand{\MR}[1]{\mbox{\href{http://www.ams.org/mathscinet-getitem?mr=#1}{#1}}}
  \renewcommand{\MR}[1]{\mbox{\href{http://www.ams.org/mathscinet-getitem?mr=#1}{#1}}}
\providecommand{\bysame}{\leavevmode\hbox to3em{\hrulefill}\thinspace}
\providecommand{\MR}{\relax\ifhmode\unskip\space\fi MR }
\providecommand{\MRhref}[2]{%
  \href{http://www.ams.org/mathscinet-getitem?mr=#1}{#2}
}
\providecommand{\href}[2]{#2}

\end{document}